\documentclass[11pt]{article}
\renewcommand{\baselinestretch}{1.2}
\usepackage{amsthm,enumerate}
\theoremstyle{plain}
\newtheorem{theorem}{Theorem}[section]
\newtheorem{lemma}[theorem]{Lemma}
\newtheorem{corollary}[theorem]{Corollary}
\newtheorem{proposition}[theorem]{Proposition}
\theoremstyle{definition}
\newtheorem{definition}[theorem]{Definition}
\theoremstyle{remark}
\newtheorem*{remark}{Remark}

\newtheorem*{rems}{Remarks}
\newtheorem*{exs}{Examples}
\newenvironment{remarks}{\begin{rems}\begin{enumerate}}{\end{enumerate}\end{rems}}

\newenvironment{items}{\begin{enumerate}[\rm (i)]}{\end{enumerate}}
\newenvironment{alphitems}{\begin{enumerate}[\rm (a)]}{\end{enumerate}}

 \usepackage{amsmath,amssymb,amsfonts,diagrams}
\newenvironment{keywords}{\noindent\small {\it Keywords\/}:}{\vskip 4pt}
\newenvironment{classification}{\noindent\small 2000 {\it Mathematics Subject
Classification\/}:}{\vskip 12pt}

\renewcommand{\implies}{\quad\Longrightarrow\quad}

\newcommand{\comps}{{\mathbb C}}

\newcommand{\posints}{{\mathbb N}}

\newcommand{\tensor}{\otimes}

\newcommand{\Tensor}{\hat{\otimes}}

\newcommand{\cstar}{{C^\ast}}

\newcommand{\id}{{\mathrm{id}}}

\newcommand{\A}{{\mathfrak A}}

\newcommand{\Hilbert}{{\mathfrak H}}

\newcommand{\varcl}[1]{\overline{#1}}

\title{Can ${\cal B}(\ell^p)$ ever be amenable?}
\author{\textit{Matthew Daws} \and \textit{Volker Runde}\thanks{Research supported by NSERC.}}
\date{}
\begin{document}
\maketitle
\begin{abstract}
It is known that ${\cal B}(\ell^p)$ is not amenable for $p =1,2,\infty$, but whether or not ${\cal B}(\ell^p)$ is amenable for $p \in (1,\infty) \setminus \{ 2 \}$ is an open problem. We show that, if ${\cal B}(\ell^p)$ is amenable for $p \in (1,\infty)$, then so are $\ell^\infty({\cal B}(\ell^p))$ and $\ell^\infty({\cal K}(\ell^p))$. Moreover, if $\ell^\infty({\cal K}(\ell^p))$ is amenable so is $\ell^\infty(\mathbb{I},{\cal K}(E))$ for any index set $\mathbb I$ and for any infinite-dimensional ${\cal L}^p$-space $E$; in particular, if $\ell^\infty({\cal K}(\ell^p))$ is amenable for $p \in (1,\infty)$, then so is $\ell^\infty({\cal K}(\ell^p \oplus \ell^2))$. We show that $\ell^\infty({\cal K}(\ell^p \oplus \ell^2))$ is not amenable for $p =1,\infty$, but also that our methods fail us if $p \in (1,\infty)$. Finally, for $p \in (1,2)$ and a free ultrafilter $\cal U$ over $\posints$, we exhibit a closed left ideal of $({\cal K}(\ell^p))_{\cal U}$ lacking a right approximate identity, but enjoying a certain, very weak complementation property.
\end{abstract}
\section*{Introduction}
In his seminal memoir \cite{Joh1}, B.\ E.\ Johnson initiated the theory of amenable Banach algebras. The choice of terminology is motivated by \cite[Theorem 2.5]{Joh1}: a locally compact group $G$ is amenable in the usual sense (see \cite{Pat}, for instance) if and only if its group algebra $L^1(G)$ is an amenable Banach algebra.
\par 
Ever since \cite{Joh1} was published, there have been ongoing efforts to determine, for particular classes of Banach algebras, which algebras in them are the amenable ones. One spectacular result in this direction is the characterization of the amenable $\cstar$-algebras: a $\cstar$-algebra is amenable if and only if it is nuclear (this result, mostly credited to A.\ Connes and U.\ Haagerup, is the culmination of the efforts of many mathematicians; see \cite{LoA} or \cite{Tak} for self-contained accounts).
\par 
One particular class of Banach algebras for which the problem of characterizing its amenable members is still wide open is the class of Banach algebras ${\cal B}(E)$, the algebras of all bounded linear operators on a Banach space $E$. From a philosophical point of view, this problem ought to be easy: amenability can often be thought of as a weak finiteness condition, and, for any infinite-dimensional Banach space $E$, the algebra ${\cal B}(E)$ should simply be too ``large'' to be amenable. Already Johnson asked in \cite{Joh1}:
\begin{itemize}
\item Is ${\cal B}(E)$ ever amenable for infinite-dimensional $E$? (\cite[10.4]{Joh1})
\item Is ${\cal B}(\Hilbert)$ amenable for an infinite-dimensional Hilbert space $\Hilbert$? (\cite[10.2]{Joh1})
\end{itemize}
\par 
The Hilbert space case was settled relatively quickly: in \cite{Was}, S.\ Wassermann showed that a nuclear von Neumann algebra had to be subhomogeneous. In view of the equivalence of amenability and nuclearity for $\cstar$-algebras, this means that ${\cal B}(\Hilbert)$ can be amenable only if $\dim \Hilbert < \infty$.
\par
Ever since, very little progress has been made in the general Banach space case. Until recently, it was not even known whether ${\cal B}(\ell^p)$ was amenable or not for any $p \in [1,\infty]$ other than $2$. This situation changed with C.\ J.\ Read's paper \cite{Rea}:
making ingenious use of random hypergraphs, Read showed that ${\cal B}(\ell^1)$ is not amenable. Moreover, he showed that, for any $p \in [1,\infty] \setminus \{ 2 \}$, the Banach algebra $\text{$\ell^\infty$-}\bigoplus_{n=1}^\infty {\cal B}(\ell^p_n)$ also fails to be amenable (the $p=2$ case already follows from Wassermann's result). Subsequently, G.\ Pisier simplified Read's proof by replacing the random hypergraphs of \cite{Rea} with expanders (\cite{Pis}). Eventually, N.\ Ozawa, simplified Pisier's argument even further and succeeded in giving a proof that simultaneously established the non-amenability of ${\cal B}(\ell^p)$ for $p=1,2,\infty$ and of $\text{$\ell^\infty$-}\bigoplus_{n=1}^\infty {\cal B}(\ell^p_n)$ for any $p \in [1,\infty]$ (\cite{Oza}); even though it is not explicitly stated in \cite{Oza}, the proof also works for ${\cal B}(c_0)$.
\par 
In the present paper, we investigate what consequences the hypothetical amenability of ${\cal B}(\ell^p)$ for $p \in (1,\infty) \setminus \{ 2 \}$ would have.
\par 
Our first result is that, if ${\cal B}(\ell^p)$ is amenable, then so is $\ell^\infty({\cal B}(\ell^p))$. As the much ``smaller'' algebra $\text{$\ell^\infty$-}\bigoplus_{n=1}^\infty {\cal B}(\ell^p_n)$ is not amenable, this lends again support to the belief that ${\cal B}(\ell^p)$ is not amenable (even though, of course, this is a far cry from a proof). 
\par
A straightforward consequence of the amenability of $\ell^\infty({\cal B}(\ell^p))$ is that $\ell^\infty({\cal K}(\ell^p))$ is amenable, too, and we shall devote most of this paper to exploring the consequences of the amenability of that particular Banach algebra and, more generally, of $\ell^\infty({\cal K}(E))$ for particular Banach spaces $E$. (Incidentally, the question of whether $\ell^\infty({\cal K}(E))$ is amenable for specific Banach spaces $E$ seems to have received almost no attention in the literature; the only references known to the authors are \cite{CSR} and \cite{LLW}, where the case $E = \ell^2$ is settled in the negative.)
\par 
First, we show that, due to the separability of ${\cal K}(\ell^p)$, the amenability of $\ell^\infty({\cal K}(\ell^p))$ already entails the amenability of $\ell^\infty(\mathbb{I}, {\cal K}(\ell^p))$ for every index set $\mathbb I$ and thus of $({\cal K}(\ell^p))_{\cal U}$ for every ultrafilter $\cal U$ (so that ${\cal K}(\ell^p)$ is ultra-amenable in the terminology of \cite{Daw}).
\par 
Next, we see that the amenability of $\ell^\infty({\cal K}(\ell^p))$ forces $\ell^\infty(\mathbb{I}, {\cal K}(E))$ to be amenable for every index set $\mathbb I$ and every infinite-dimensional ${\cal L}^p$-space $E$ in the sense of \cite{LP}. In particular, if $\ell^\infty({\cal K}(\ell^p))$ is amenable, the so is $\ell^\infty({\cal K}(\ell^p \oplus \ell^2))$, which is interesting because ${\cal B}(\ell^p \oplus \ell^2)$ is known to be non-amenable.
\par 
We then study the amenability of $\ell^\infty({\cal K}(E \oplus F))$ for certain Banach spaces $E$ and $F$. Using the theory of operator ideals (see \cite{Pie}), we show that $\ell^\infty({\cal K}(E \oplus \ell^2))$ is not amenable for $E = c_0, \ell^\infty$, and $\ell^1$, but we also show that our methods fail to establish the non-amenability of $\ell^\infty({\cal B}(\ell^p \oplus \ell^2))$ for $p \in (1,\infty)$.
\par 
Finally, we take a look at a particular left ideal of $({\cal K}(\ell^p))_{\cal U}$ for $p \in (1,2)$ and $\cal U$ a free ultrafilter over $\posints$. We show that this ideal lacks a right approximate identity and, at the same time, enjoys a certain complementation property, which is unfortunately too weak to obtain a contradiction to the amenability of $({\cal K}(\ell^p))_{\cal U}$.
\subsubsection*{Acknowledgments}
This research was initiated while the first author was visiting the University of Alberta in the summer of 2007; the financial support and the kind hospitality are gratefully acknowledged. Both authors would like to thank Andreas Defant, Albrecht Pietsch, and Nicole Tomczak-Jaegermann for valuable help with operator ideals.
\section{Amenable Banach algebras} 
The definition of an amenable Banach algebra given in \cite{Joh1} is in terms of certain derivations being inner. Throughout this paper, however, we shall not rely on that definition directly, but rather on a more intrinsic, but equivalent characterization, also due to Johnson (\cite{Joh2}).
\par
Let $\A$ be a Banach algebra, and let $E$ and $F$ be a left and right Banach $\A$-module, respectively. We use $\Tensor$ to denote the projective tensor product of Banach spaces. The Banach space $E \Tensor F$ becomes a Banach $\A$-bimodule via
\[
  a \cdot ( x \tensor y) := a \cdot x \tensor y \quad\text{and}\quad (x \tensor y) \cdot a : =x \tensor y \cdot a
  \qquad (a \in \A, \, x \in E, \, y \in F).
\]
In particular, $\A \Tensor \A$ is a Banach $\A$-bimodule in a canonical manner. With respect to these module operatations, the diagonal map $\Delta \!: \A \Tensor \A \to \A$ induced by multiplication, i.e., $\Delta(a \tensor b) = ab$ for $a,b \in \A$, is a bimodule homomorphism; if we want to emphasize the algebra $\A$, we sometimes write $\Delta_\A$ for $\Delta$.
\begin{definition} \label{amdef}
Let $\A$ be a Banach algebra. An \emph{approximate diagonal} for $\A$ is a bounded net $( \boldsymbol{d}_\alpha )_\alpha$ in $\A \Tensor \A$ such that
\begin{equation} \label{diag1}
  a \cdot \boldsymbol{d}_\alpha - \boldsymbol{d}_\alpha \cdot a \to 0 \qquad (a \in \A)
\end{equation}
and
\begin{equation} \label{diag2}
  a \Delta \boldsymbol{d}_\alpha \to a \qquad (a \in \A).
\end{equation}
If $\A$ has an approximate diagonal bounded, we say that $\A$ is \emph{amenable}. 
\end{definition}
\begin{remarks}
\item If $\A$ is amenable and has an identity $1_\A$, then there is an approximate diagonal $( \boldsymbol{d}_\alpha )_{\alpha \in \mathbb A}$ for $\A$ such that $\Delta \boldsymbol{d}_\alpha = 1_\A$ for all $\alpha \in \mathbb A$.
\item If $\A$ is amenable with an approximate diagonal bounded by $C \geq 1$, then $\A$ is also called \emph{$C$-amenable}. It is clear from (\ref{diag2}) that is doesn't make sense to speak of $C$-amenability for any $C \in (0,1)$.
\end{remarks}
\par 
For modern accounts of the theory of amenable Banach algebras, see \cite{Dal} or \cite{LoA}. 
\par 
The hereditary properties of Banach algebraic amenability are well understood. For instance, quotients of amenable Banach algebras are again amenable (\cite[Corollary 2.3.2]{LoA}), and a closed ideal of an amenable Banach algebra is amenable if and only if it has a bounded approximate identity and if and only if it is weakly complemented (\cite[Theorem 2.3.7]{LoA}). Whether or not a particular subalgebra of an amenable Banach algebra is amenable is a much more delicate question and an elegant characterization is certainly out of reach. Nevertheless, some partial results exist such as the following (\cite[Theorem 6.2]{GJW}):
\begin{theorem} \label{bgn}
Let $\A$ be a Banach algebra with a bounded approximate identity, let $P_1 \in {\cal M}(\A)$ be a projection, and let $P_2 := \id_\A - P_1$. Suppose further that $\Delta_\A$ maps $P_2 \A P_1 \Tensor P_1 \A P_2$ onto $P_2 \A P_2$. Then $\A$ is amenable if and only if $P_1 \A P_1$ is amenable.
\end{theorem}
\par 
Here, ${\cal  M}(\A)$ stands for the \emph{multiplier algebra} of $\A$ (\cite[p.\ 60]{Dal}).
\par 
The Banach algebra $\A$ in Theorem \ref{bgn} has a matrix like structure thanks to the projections $P_1$ and $P_2$. We shall now prove a necessary condition for the non-amenability of such algebras:
\begin{proposition} \label{noamprop}
Let $\A$ be a Banach algebra, let $P_1 \in {\cal M}(\A)$ be an idempotent, and let $P_2 := \id_\A - P_1$. Suppose that there is a closed ideal $I$ of $\A$ such that $P_2 I P_1 = P_2 \A P_1$, but $P_1 I P_2 \subsetneq P_1 \A P_2$. Then $\A$ is not amenable.
\end{proposition}
\begin{proof}
For $j,k \in \{1,2\}$, set
\[
  \A_{j,k} := P_j \A P_k \qquad\text{and}\qquad I_{j,k} := P_j I P_k. 
\]
It follows that
\[
  \A \cong \begin{bmatrix} \A_{1,1} & \A_{1,2} \\ \A_{2,1} & \A_{2,2} \end{bmatrix}
  \qquad\text{and}\qquad
  I \cong \begin{bmatrix} I_{1,1} & I_{1,2} \\ I_{2,1} & I_{2,2} \end{bmatrix}
\]
and thus
\[
  \A / I \cong \begin{bmatrix} \A_{1,1} / I_{1,1} & \A_{1,2} / I_{1,2} \\ \A_{2,1} / I_{2,1} & \A_{2,2} / I_{2,2} \end{bmatrix}
  = \begin{bmatrix} \A_{1,1} / I_{1,1} & \A_{1,2} / I_{1,2} \\ 0 & \A_{2,2} / I_{2,2} \end{bmatrix}.
\]
Consequently, $\left[ \begin{smallmatrix} 0 & \A_{1,2} / I_{1,2} \\ 0 & 0 \end{smallmatrix} \right]$ is a non-zero, complemented, nilpotent ideal of $\A/I$, which is impossible if $\A/I$ is amenable.
\end{proof}
\begin{remark}
The idea behind Proposition \ref{noamprop} is implicitly already contained in \cite[Question 4]{Gro}, where it is attributed to G.\ A.\ Willis. It can be used to establish the non-amenability of ${\cal B}(E \oplus F)$ if ${\cal B}(F,E) = {\cal K}(F,E)$, but ${\cal B}(E,F) \subsetneq {\cal K}(E,F)$; this applies, for instance, to ${\cal B}(\ell^p \oplus \ell^q)$ with $1 \leq p < q < \infty$ (\cite[Proposition 2.c.3]{LT}).
\end{remark}
\par 
We conclude this section with the discussion of a stronger variant of amenability also introduced by Johnson (\cite{Joh3}).
\par 
Given a Banach algebra $\A$, let $\Sigma$ denote the flip map on $\A \Tensor \A$, i.e., $\Sigma(a \tensor b) = b \tensor a$ for $a,b \in \A$. An element $\boldsymbol{a} \in \A \tensor \A$ is called \emph{symmetric} if $\Sigma \boldsymbol{a} = \boldsymbol{a}$; somewhat abusing terminology, we will also call a net in $\A \tensor \A$ symmetric if it consists of symmetric elements of $\A \Tensor \A$.
\begin{definition}
A Banach algebra is called \emph{symmetrically amenable} if it has a symmetric approximate diagonal.
\end{definition}
\begin{remark}
The group algebra $L^1(G)$ of a locally compact group $G$ is symmetrically amenable if and only if it is amenable (\cite[Theorem 4.1]{Joh3}) whereas the Cuntz algebras ${\cal O}_n$ for $n \in \posints$, $n \geq 2$ (\cite{Cun}) are amenable, but not symmetrically amenable (\cite[p.\ 457]{Joh3}).
\end{remark}
\par 
In view of how difficult it is, even for very well behaved Banach spaces $E$, to show that ${\cal B}(E)$ is not amenable, it is somewhat surprising to see how easily the correspoding question for symmetric amenability can be settled in the negative for a large class of Banach spaces:
\begin{proposition}
Let $E$ be a Banach space such that $E \cong E \oplus E$. Then ${\cal B}(E)$ is not symmetrically amenable.
\end{proposition}
\begin{proof}
Assume that ${\cal B}(E)$ is symmetrically amenable. By \cite[Corollary 2.5]{Joh3}, there is $\phi \in {\cal B}(E)^\ast$ such that
$\langle \id_E, \phi \rangle = 1$ and $\langle ST, \phi \rangle = \langle TS, \phi \rangle$ for $S,T \in {\cal B}(E)$.
\par 
For $j=1,2$, let $P_j \!: E \oplus E \to E$ be the projection onto the $j$-th summand. Since $E \cong E \oplus E$, there are $U_j,V_j \in {\cal B}(E)$ with
\[
  U_j V_j = \id_E \quad\text{and}\quad V_j U_j = P_j \qquad (j=1,2).
\]
It follows that
\[
  1 = \langle \id_E,\phi \rangle = \langle P_1 + P_2 , \phi \rangle 
    = \langle V_1 U_1, \phi \rangle + \langle V_2 U_2, \phi \rangle 
    = \langle U_1 V_1 , \phi \rangle + \langle U_2 V_2 , \phi \rangle
    = 2 \langle \id_E,\phi \rangle = 2,
\]
which is nonsense.
\end{proof}
\section{Amenability of ${\cal B}(\ell^p)$ and $\ell^\infty({\cal K}(\ell^p))$}
We begin with establishing some notation, part of which was already used in the introduction.
\par 
Let $\mathbb I$ be any index set, and let $( E_i )_{i \in \mathbb I}$ be a family of Banach spaces; we write $\prod_{i \in \mathbb I} E_i$ for its Cartesian product. For $p \in [1,\infty)$, we set
\[
  \text{$\ell^p$-}\bigoplus_{i \in \mathbb I} E_i 
  := \left\{ (x_i)_{i \in \mathbb I} \in \prod_{i \in \mathbb I} E_i : \sum_{i \in \mathbb I} \| x_i \|^p < \infty \right\};
\]
it is a linear space which becomes a Banach space if equipped with the norm
\[
  \| (x_i)_{i \in \mathbb I} \|_p := \left( \sum_{i \in \mathbb I} \| x_i \|^p \right)^\frac{1}{p}
  \qquad \left( (x_i)_{i \in \mathbb I} \in \text{$\ell^p$-}\bigoplus_{i \in \mathbb I} E_i \right).
\]
Furthermore, we define
\[
  \text{$\ell^\infty$-}\bigoplus_{i \in \mathbb I} E_i := 
  \left\{ (x_i)_{i \in \mathbb I} \in \prod_{i \in \mathbb I} E_i : \sup_{i \in \mathbb I} \| x_i \| < \infty \right\};
\]
it, too, becomes a Banach space with the norm
\[
  \| (x_i)_{i \in \mathbb I} \|_\infty := \sup_{i \in \mathbb I} \| x_i \|  \qquad 
  \left( (x_i)_{i \in \mathbb I} \in \text{$\ell^\infty$-}\bigoplus_{i \in \mathbb I} E_i \right).
\]
We use $\text{$c_0$-}\bigoplus_{i \in \mathbb I} E_i$ to denote the closure of those $( x_i )_{i \in \mathbb I} \in \text{$\ell^\infty$-}\bigoplus_{i \in \mathbb I} E_i$ for which $x_i = 0$ for all but finitely many $i \in \mathbb I$. We note that, if $( \A_i )_{i \in \mathbb I}$ is a family of Banach algebras, then $\text{$\ell^\infty$-}\bigoplus_{i \in \mathbb I} \A_i$ is a Banach algebra (which contains $\text{$c_0$-}\bigoplus_{i \in \mathbb I} \A_i$ as a closed ideal). If $E_i = E$ for all $i \in \mathbb I$, we simply write $\ell^p(\mathbb{I},E)$ or $c_0(\mathbb{I},E)$ instead of $\text{$\ell^p$-}\bigoplus_{i \in \mathbb I}  E$ and $\text{$c_0$-}\bigoplus_{i \in \mathbb I} E$, respectively. We apply the usual conventions: if $E_i = \comps$ for all $i \in \mathbb I$ or $\mathbb{I} = \posints$, we suppress the symbol for the space or the index set, respectively. For instance, if $p \in [1,\infty]$ and $E$ is any Banach space, then $\ell^p(E)$ stands for $\ell^p(\posints,E)$, and if $\mathbb{I}$ is any index set, then  $c_0(\mathbb{I})$ means $c_0(\mathbb{I},\comps)$. Also, we write $\ell^p_n$ instead of $\ell^p(\{ 1, \ldots, n \},\comps)$. Finally, for $i \in \mathbb I$, we let $\delta_i \!: \mathbb{I} \to \comps$ denote the point mass at $i$; it is clear that $\delta_i \in \ell^p(\mathbb{I})$ for any $p \in [1,\infty]$.
\par 
Given any Banach space $E$, we have isometric isomorphisms between $\ell^p(\ell^p(E)) = \ell^p(\posints^2,E)$ and $\ell^p(E)$ for $p \in [1,\infty)$ and between $c_0(c_0(E)) = c_0(\posints^2,E)$ and $c_0(E)$ (simply due to the fact that $\posints$ and $\posints^2$ have the same cardinality). This simple observation lies at the heart of the proof of our first theorem:
\begin{theorem} \label{thm1}
Let $E$ be a Banach space. Then:
\begin{items}
\item for $p \in [1,\infty)$, the Banach algebra ${\cal B}(\ell^p(E))$ is amenable if and only if $\ell^\infty({\cal B}(\ell^p(E)))$ is amenable;
\item ${\cal B}(c_0(E))$ is amenable if and only if $\ell^\infty({\cal B}(c_0(E)))$ is amenable.
\end{items}
\end{theorem}
\begin{proof}
We only prove (i) ((ii) is proven analogously).
\par 
Set $\A := \ell^\infty({\cal B}(\ell^p(E)))$. It is elementary that ${\cal B}(\ell^p(E))$ is amenable if $\A$ is (\cite[Corollary 2.3.2]{LoA}). 
\par 
For the converse, suppose that ${\cal B}(\ell^p(E))$ is amenable, and let $( \boldsymbol{d}_\alpha )_{\alpha \in \mathbb A}$ be an approximate diagonal for it; we may suppose that $\Delta \boldsymbol{d}_\alpha = \id_{\ell^p(E)}$ for all $\alpha \in \mathbb A$.
\par 
First, observe that we can identify $\A$ with the block diagonal matrices in ${\cal B}(\ell^p(\ell^p(E)))$. For $n \in \posints$, let $P_n \!: \ell^p(\ell^p(E)) \to \ell^p(E)$ denote the projection onto the $n$-th coordinate. Define
\[
  {\cal Q} \!: {\cal B}(\ell^p(\ell^p(E))) \to {\cal B}(\ell^p(\ell^p(E))), \quad T \mapsto \sum_{n=1}^\infty P_n T P_n,
\] 
where the infinite series converges in the strong operator topology $\mathrm{SOT}$. Then $\cal Q$ is a projection onto $\A$.
\par 
Since we have $\ell^p(\ell^p(E)) \cong \ell^p(E)$, there are bounded sequences $(U_n)_{n=1}^\infty$ and $( V_n )_{n=1}^\infty$ in ${\cal B}(\ell^p(\ell^p(E)))$ such that
\[
  U_n V_n = \id_{\ell^p(\ell^p(E))} \quad\text{and}\quad V_n U_n = P_n \qquad (n \in \posints)
\]
and, consequently, 
\begin{equation} \label{UV}
  V_n = P_n V_n \quad\text{and}\quad U_n = U_n P_n \qquad (n \in \posints).
\end{equation}
Define
\begin{align*}
  {\cal Q}_L \!: {\cal B}(\ell^p(\ell^p(E))) \to {\cal B}(\ell^p(\ell^p(E))), & \quad T \mapsto \sum_{n=1}^\infty P_n T U_n \\
  \intertext{and}
  {\cal Q}_R \!: {\cal B}(\ell^p(\ell^p(E))) \to {\cal B}(\ell^p(\ell^p(E))), & \quad T \mapsto \sum_{n=1}^\infty V_n T P_n,
\end{align*}
where again the series are convergent in the strong operator topology. It is obvious that ${\cal Q}_L$ is a left and ${\cal Q}_R$ a right $\A$-module homomorphism, and with (\ref{UV}) in mind, it is easy to see that both ${\cal Q}_L$ and ${\cal Q}_R$ attain their values in $\A$.
\par 
Let $S,T \in {\cal B}(\ell^p(\ell^p(E)))$. Since multiplication is jointly continuous on norm bounded subsets with respect to $\mathrm{SOT}$, we have
\[
  \begin{split}
  ({\cal Q}_L S) ({\cal Q}_R T) & = \text{$\mathrm{SOT}$-}\lim_{N\to\infty} 
                 \left( \sum_{n=1}^N P_n S U_n \right)\left( \sum_{n=1}^N V_n T P_n \right) \\
  & = \text{$\mathrm{SOT}$-}\lim_{N\to\infty} \sum_{n,m=1}^N P_n S U_n V_m T P_m \\
  & = \text{$\mathrm{SOT}$-}\lim_{N\to\infty} \sum_{n=1}^N P_n S U_n V_n T P_n \\
  & = \text{$\mathrm{SOT}$-}\lim_{N\to\infty} \sum_{n=1}^N P_n ST P_n \\
  & = {\cal Q}(ST). 
  \end{split}
\]
It follows that
\begin{equation} \label{Q}
  \Delta \circ ({\cal Q}_L \tensor {\cal Q}_R) = {\cal Q} \circ \Delta.
\end{equation}
\par 
Identifying ${\cal B}(\ell^p(E))$ and ${\cal B}(\ell^p(\ell^p(E))$, we claim that $( ({\cal Q}_L \tensor {\cal Q}_R) \boldsymbol{d}_\alpha )_{\alpha \in \mathbb A}$ is an approximate diagonal for $\A$. Clearly, the net is bounded in $\A \Tensor \A$. Since ${\cal Q}_L \tensor {\cal Q}_R$ is an $\A$-bimodule homomorphism, it follows that (\ref{diag1}) holds. Finally, since
\[
  \Delta(({\cal Q}_L \tensor {\cal Q}_R) \boldsymbol{d}_\alpha) = {\cal Q}(\Delta \boldsymbol{d}_\alpha ) = 
  {\cal Q}\left(\id_{\ell^p(\ell^p(E))}\right) = \id_{\ell^p(\ell^p(E))}
\]
by (\ref{Q}), condition (\ref{diag2}) holds as well.
\end{proof}
\par 
Specializing to $E = \comps$ yields:
\begin{corollary} \label{cor1}
Let $p \in (1,\infty)$ be such that ${\cal B}(\ell^p)$ is amenable. Then $\ell^\infty({\cal B}(\ell^p))$ and $\ell^\infty({\cal K}(\ell^p))$ are both amenable.
\end{corollary}
\begin{proof}
The claim for $\ell^\infty({\cal B}(\ell^p))$ is an immediate consequence of Theorem \ref{thm1}. Since ${\cal K}(\ell^p)$ has a bounded approximate identity, so does $\ell^\infty({\cal K}(\ell^p))$. Since $\ell^\infty({\cal K}(\ell^p))$ is a closed ideal of the amenable Banach algebra $\ell^\infty({\cal B}(\ell^p))$, it is amenable by \cite[Proposition 2.3.3]{LoA}.
\end{proof}
\begin{remarks}
\item The analogous statement of Corollary \ref{cor1} for $\ell^1$ and $c_0$ is also true. Since, however, ${\cal B}(\ell^1)$ and ${\cal B}(c_0)$ are known to be not amenable by \cite{Oza}, it would be somewhat pointless to formulate it.
\item Even though ${\cal B}(\ell^1)$ and ${\cal B}(c_0)$ are not amenable, it seems to be unknown whether $\ell^\infty({\cal K}(\ell^1))$ and $\ell^\infty({\cal K}(c_0))$ are.
\item It is known that $\ell^\infty({\cal K}(\ell^2))$ is not amenable (see \cite{LLW}), but proving it is at about the same level of difficulty as a proof for the non-amenability of ${\cal B}(\ell^2)$.
\end{remarks}
\par 
We shall thus, from now on, focus on the (non-)amenability of $\ell^\infty({\cal K}(\ell^p))$ instead of that of ${\cal B}(\ell^p)$.
\section{Ultra-amenability of ${\cal K}(\ell^p)$}
Let $E$ be a Banach space, let $\mathbb I$ be an index set, and let $\cal U$ be an ultrafilter over $\mathbb I$. We let
\[
  {\cal N}_{\cal U} := \left\{ (x_i)_{i \in \mathbb I} : \lim_{i \in \cal U} \| x_i \| = 0 \right\}.
\] 
It is immediate that ${\cal N}_{\cal U}$ is a closed subspace of $\ell^\infty(\mathbb{I},E)$. The quotient space $\ell^\infty(\mathbb{I},E)/ {\cal N}_{\cal U}$ is called the \emph{ultrapower of $E$ with respect to $\cal U$}; we denote it by $(E)_{\cal U}$. Whenever $(x_i)_{i \in \mathbb I} \in \ell^\infty(\mathbb{I},E)$, we write $( x_i )_{\cal U}$ for its equivalence class in $\cal U$. For further material on ultrapowers, we refer to the survey article \cite{Hei} and the somewhat more detailed treatment in \cite{Sim}. 
\par 
If $\A$ is a Banach algebra, then it is straightforward that $(\A)_{\cal U}$ is again a Banach algebra. The following definition is due to the first author (\cite{Daw}):
\begin{definition} \label{uldef}
A Banach algebra $\A$ is said to be \emph{ultra-amenable} if $(\A)_{\cal U}$ is amenable for every ultrafilter $\cal U$. 
\end{definition}
\begin{remark}
Ultra-amenability implies amenability (\cite[Corollary 5.5]{Daw}), but is much stron\-ger: a $\cstar$-algebra is ultra-amenable if and  only if it is subhomogeneous (\cite[Theorem 5.7]{Daw}) and $\ell^1(G)$, for a discrete group $G$, is ultra-amenable if and only if $G$ is finite (\cite[Theorem 5.11]{Daw}).
\end{remark}
\par 
Suppose that ${\cal B}(\ell^p)$ is amenable for some $p \in (1,\infty)$. Then $\ell^\infty({\cal K}(\ell^p))$ is amenable by Corollary \ref{cor1}, so that its quotient $({\cal K}(\ell^p))_{\cal U}$ is amenable for every ultrafilter $\cal U$ over $\posints$. Alas, this does not allow us (yet) to say that ${\cal K}(\ell^p)$ is ultra-amenable because Definition \ref{uldef} requires us to consider ultrafilters over arbitrary index sets. 
\par 
Nevertheless, the amenability of $\ell^\infty({\cal K}(\ell^p))$ allows us to conclude the ultra-amenability of ${\cal K}(\ell^p)$ by virtue of the following theorem:
\begin{theorem} \label{thm2}
The following are equivalent for a \emph{separable} Banach algebra $\A$:
\begin{items}
\item $\ell^\infty(\A)$ is amenable;
\item $\ell^\infty(\mathbb{I},\A)$ is amenable for every index set $\mathbb I$;
\item $\A$ is ultra-amenable.
\end{items}
\end{theorem}
\par 
For the proof, recall the following definitions from \cite{Daw}. Let $\A$ be a Banach algebra, and let $n \in \posints$. Then:
\begin{itemize}
\item let $S_n(\A)$ denote the collection of all subsets of the unit sphere of $\A$ of cardinality $n$;
\item for $C \geq 1$ and $\epsilon > 0$, let $D_n(\A,C,\epsilon)$ consist of those $A \subset S_n(\A)$ such that there is a sequence $( t_k )_{k=1}^\infty$ in $[0,\infty)$ with $\sum_{k=1}^\infty t_k \leq C$ with the property that, for each $S \in A$, there are sequences $( a_k )_{k=1}^\infty$ and $( b_k )_{k=1}^\infty$ in $\A$ with $\| a_k \| \| b_k \| \leq t_k$ for $k \in \posints$, so that
\[
  \boldsymbol{d} := \sum_{k=1}^\infty a_k \tensor b_k \in \A \Tensor \A,
\]
and
\[
  \| a \cdot \boldsymbol{d} - \boldsymbol{d} \cdot a \| \leq \epsilon 
  \quad\text{and}\quad
  \| \Delta_\A(\boldsymbol{d}) a - a \| \leq \epsilon \qquad (a \in F).
\]
\end{itemize}
\begin{lemma} \label{ullem}
For a Banach algebra $\A$ consider the following statements:
\begin{items}
\item there is $C \geq 1$ such that $S_n(\A) \in D_n(\A,C,\epsilon)$ for each $n \in \posints$ and $\epsilon > 0$;
\item $\ell^\infty(\mathbb{I},\A)$ is amenable for each index set $\mathbb I$;
\item $\ell^\infty(\A)$ is amenable.
\end{items}
Then 
\[
  \emph{(i)} \implies \emph{(ii)} \implies \emph{(iii)}, 
\]
and \emph{(iii)} implies \emph{(i)} if $\A$ is separable.
\end{lemma}
\begin{proof}
(i) $\Longrightarrow$ (ii) is routine in view of the definition of $D_n(\A,C,\epsilon)$, and (ii) $\Longrightarrow$ (iii) is trivial.
\par 
Suppose that (iii) holds and that $\A$ is separable. Let $C \geq 1$ be such that $\A$ is $C$-amenable, let $n \in \posints$, and let $\epsilon > 0$. 
\par 
Define a metric $d$ on $S_n(\A)$ by letting 
\[
  d(A,B) := \max_{a \in A} \min_{b \in B} \| a - b \| + \max_{b \in B} \min_{a \in A} \| a - b \|\qquad (A,B \in S_n(\A)).
\]
The separability of $\A$ implies that the metric space $(S_n(\A),d)$ is separable and so contains a dense, countable subset, say $\{ A_1, A_2, \ldots \}$. For each $k \in \posints$, let $A_k = \left\{ a^{(k)}_1, \ldots, a^{(k)}_n \right\}$. For $j =1, \ldots, n$, set $a_j := \left( a^{(k)}_j \right)_{k=1}^\infty \in \ell^\infty(\A)$. Let $( b_k )_{k=1}^\infty$ and $( c_k )_{k=1}^\infty$ be sequences in $\ell^\infty(\A)$ with $\sum_{k=1}^\infty \| b_k \| \| c_k \| \leq C$ such that, for $\boldsymbol{d} := \sum_{k=1}^\infty b_k \tensor c_k \in \ell^\infty(\A) \Tensor \ell^\infty(\A)$, we have
\begin{equation} \label{diag3}
  \| a_j \cdot \boldsymbol{d} - \boldsymbol{d} \cdot a_j \| \leq \frac{\epsilon}{2}
\end{equation}
and
\begin{equation} \label{diag4}
  \| \Delta_{\ell^\infty(\A)} ( \boldsymbol{d}) a_j - a_j \| \leq \frac{\epsilon}{2}
\end{equation}
for $j =1, \ldots, n$. (The existence of such $( b_k )_{k=1}^\infty$ and $( c_k )_{k=1}^\infty$ follows from the $C$-amenability of $\ell^\infty(\A)$.) 
\par 
For each $k \in \posints$, let $b_k = \left( b^{(k)}_\nu \right)_{\nu=1}^\infty$ and $c_k = \left( c^{(k)}_\nu \right)_{\nu=1}^\infty$. Then (\ref{diag4}) yields that
\begin{equation} \label{diag5}
  \sup_{\nu \in \posints} \left\| \sum_{k=1}^\infty b^{(k)}_\nu c^{(k)}_\nu a_j^{(\nu)}  - a_j^{(\nu)} \right\| \leq \frac{\epsilon}{2}
  \qquad (j = 1, \ldots, n).
\end{equation}
For $\nu \in \posints$, let $P_\nu \!: \ell^\infty(\A) \to \A$ be the projection onto the $\nu$-th coordinate, and note that $P_\nu \tensor P_\nu \!: \ell^\infty(\A) \Tensor \ell^\infty(\A) \to \A \Tensor \A$ is a contraction. From (\ref{diag3}), we conclude that
\begin{equation} \label{diag6}
  \sup_{\nu \in \posints} \left\| \sum_{k=1}^\infty a_j^{(\nu)} b^{(k)}_\nu \tensor c^{(k)}_\nu -
     \sum_{k=1}^\infty b^{(k)}_\nu \tensor c^{(k)}_\nu a_j^{(\nu)} \right\| \leq \frac{\epsilon}{2}
  \qquad (j = 1, \ldots, n).
\end{equation}
Finally, it is straightforward that
\[
  \sum_{k=1}^\infty \sup_{\nu \in \posints} \left\| b_\nu^{(k)} \right\|\left\| b_\nu^{(k)} \right\| \leq
  \sum_{k=1}^\infty \| b_k \| \| c_k \| \leq C.
\]
\par 
Let $A \in S_n(\A)$ be arbitrary. Since $\{ A_1, A_2, \ldots \}$ is dense in $S_n(\A)$, there is $\nu \in \posints$ such that $d(A,A_\nu) \leq \frac{\epsilon}{4C}$. With $A_\nu = \left\{ a^{(\nu)}_1, \ldots, a^{(\nu)}_n \right\}$, this means that, for any $a \in A$, there is $j \in \{1, \ldots, n\}$ such that $\left\| a_j^{(\nu)} - a \right\| \leq \frac{\epsilon}{4C}$. From (\ref{diag6}), we infer that
\[
  \begin{split}
  \lefteqn{\left\| \sum_{k=1}^\infty a b^{(k)}_\nu \tensor c^{(k)}_\nu - \sum_{k=1}^\infty b^{(k)}_\nu \tensor c^{(k)}_\nu a \right\|} 
  & \\
  & \leq 2 \left\| a - a_j^{(\nu)} \right\| \sum_{k=1}^\infty \left\| b_\nu^{(k)} \right\|\left\| b_\nu^{(k)} \right\| +
  \left\| \sum_{k=1}^\infty a_j^{(\nu)} b^{(k)}_\nu \tensor c^{(k)}_\nu -
     \sum_{k=1}^\infty b^{(k)}_\nu \tensor c^{(k)}_\nu a_j^{(\nu)} \right\| \\
  & \leq \frac{\epsilon}{4C} 2C + \frac{\epsilon}{2} \\
  & = \epsilon
  \end{split}
\]
and from (\ref{diag5}) that
\[
  \begin{split}
  \lefteqn{\left\| \sum_{k=1}^\infty b^{(k)}_\nu c^{(k)}_\nu a  - a \right\|} & \\
  & \leq \left\| a - a_j^{(\nu)} \right\| + 
       \left\| a - a_j^{(\nu)} \right\|\sum_{k=1}^\infty \left\| b_\nu^{(k)} \right\|\left\| b_\nu^{(k)} \right\| + 
       \left\| \sum_{k=1}^\infty b^{(k)}_\nu c^{(k)}_\nu a_j^{(\nu)}  - a_j^{(\nu)} \right\| \\
  & \leq \frac{\epsilon}{4C} (1+C) + \frac{\epsilon}{2} \\
  & = \epsilon.
  \end{split}
\]
\par 
All in all, we have established that $S_n(\A) \in D_n(\A,C,\epsilon)$.
\end{proof}
\begin{proof}[Proof of Theorem \emph{\ref{thm2}}]
By Lemma \ref{ullem}, (i) $\Longleftrightarrow$ (ii) holds, and (ii) $\Longrightarrow$ (iii) is trivial.
\par 
We postpone the actual proof of (iii) $\Longrightarrow$ (i) for some preliminary considerations.
\par 
Let $\cal S$ be the set of all sequences in $[0,\infty)$. For $( t_k )_{k=1}^\infty, ( s_k )_{k=1}^\infty \in \cal S$, define $( t_k )_{k=1}^\infty \ll ( s_k )_{k=1}^\infty$ if there are $( r_k )_{k=1}^\infty \in \cal S$, a bijection $\sigma \!: \posints \to \posints$, and $\nu,\mu \in \posints$ such that
\begin{itemize}
\item $s_k = r_{\sigma(k)}$ for $k \in \posints$,
\item $r_k = t_k$ for $k < \nu$,
\item $r_\nu + r_{\nu+1} + \cdots + r_{\nu+\mu-1} = t_\nu$, and
\item $r_{k+\mu-1} = t_k$ for $k > \nu$.
\end{itemize}
(Informally, one might want to say that $( r_k )_{k=1}^\infty$ is obtained from $( t_k )_{k=1}^\infty$ by splitting up one term and from $( s_k )_{k=1}^\infty$ through rearrangement.) It is clear from this definition that, whenever $( t_k )_{k=1}^\infty \ll ( s_k )_{k=1}^\infty$ and one of the sequences lies in $\ell^1$, the so does the other and has the same norm in $\ell^1$. We then define $( t_k )_{k=1}^\infty \preceq ( s_k )_{k=1}^\infty$ if there are $\left( r_k^{(1)} \right)_{k=1}^\infty, \ldots, \left( r_k^{(m)} \right)_{k=1}^\infty \in \cal S$ such that
\[
  ( t_k )_{k=1}^\infty \ll \left( r_k^{(1)} \right)_{k=1}^\infty \ll \cdots \ll
  \left( r_k^{(m)} \right)_{k=1}^\infty \ll ( s_k )_{k=1}^\infty.
\]
Let ${\cal S}_0$ be the collection of all sequences in $\cal S$ that are eventually zero, and note that $({\cal S}_0, \preceq)$ is a directed set.
\par 
Let $n \in \posints$, let $F \in S_n(\A)$, let $\epsilon > 0$, and let $( t_k )_{k=1}^\infty \in \cal \ell^1 \cap S$. We call $F$ and $( t_k )_{k=1}^\infty$ \emph{compatible} if there is $\boldsymbol{d} = \sum_{k=1}^\infty a_k \tensor b_k \in \A \Tensor \A$ with $\| a_k \| \| b_k \| \leq t_k$ for $k \in \posints$ and
\[
  \| a \cdot \boldsymbol{d} - \boldsymbol{d} \cdot a \| \leq \epsilon \quad\text{and}\quad
  \| \Delta(\boldsymbol{d}) a - a \| \qquad (a \in F).
\]
Note that, if $( t_k )_{k=1}^\infty \preceq ( s_k )_{k=1}^\infty$ and $F$ and $( t_k )_{k=1}^\infty$ are compatible, then so are $F$ and $( s_k )_{k=1}^\infty$.
\par 
Suppose now that $\A$ is ultra-amenable, let $C \geq 1$ be as in \cite[Theorem 5.6]{Daw}, let $n \in \posints$, and let $\epsilon > 0$. By \cite[Theorem 5.6]{Daw}, there is a partition of $S_n(\A)$ into finitely many sets each of which has a compatible sequence in $\ell^1  \cap \cal S$ with $\ell^1$-norm at most $C$. We can suppose that each of these sequences lies in ${\cal S}_0$ and use the directedness of $({\cal S}_0, \preceq)$ to obtain one single sequence in ${\cal S}_0$---still with $\ell^1$-norm at most $C$---that is compatible with all sets in the partition of $S_n(\A)$. But then $S_n(\A)$ itself and that sequence are compatible, which means that $S_n(\A) \in D_n(\A,C,\epsilon)$. By Lemma \ref{ullem}, this implies the amenability of $\ell^\infty(\A)$.
\end{proof}
\begin{remark}
The proof of Theorem \ref{thm2} can be modified to yield the following generalization: a Banach algebra $\A$ of density character $\kappa$ is ultra-amenable if and only if $\ell^\infty(\mathbb{I},\A)$ is amenable for any index set $\mathbb I$ and if and only if it is amenable for an index set $\mathbb I$ of cardinality $\kappa$. This closes a gap in the proof of \cite[Theorem 2.5]{LLW}, which claims the equivalence of Theorem \ref{thm2}(ii) and (iii) in the case of a $\cstar$-algebra.
\end{remark}
\par 
As ${\cal K}(\ell^p)$ for $p \in (1,\infty)$ and ${\cal K}(c_0)$ are separable, Theorem \ref{thm2} yields:
\begin{corollary} \label{cor2}
Let $E = \ell^p$ with $p \in (1,\infty)$ or $E = c_0$. Then the following are equivalent:
\begin{items}
\item $\ell^\infty({\cal K}(E))$ is amenable;
\item $\ell^\infty(\mathbb{I},{\cal K}(E))$ is amenable for every index set $\mathbb I$;
\item ${\cal K}(E)$ is ultra-amenable.
\end{items}
\end{corollary}
\section{Amenability of $\ell^\infty({\cal K}(E))$ for ${\cal L}^p$-spaces}
Whether or not the Banach algebras of the form $\ell^\infty({\cal K}(E))$ for a Banach space $E$ are amenable or not seems to have received very little attention in the literature so far. As Corollary \ref{cor1} shows, it is, for $E = \ell^p$ with $p \in [1,\infty)$, intimately linked to the open problem of whether ${\cal B}(\ell^p)$ is amenable and thus certainly a question deserving further exploration.
\par 
In this section, we show that the (possible) amenability of $\ell^\infty({\cal K}(\ell^p))$ entails the amenability of $\ell^\infty({\cal K}(E))$ for a large class of a Banach spaces $E$.
\par 
For our first proposition, we denote by ${\cal F}(E,F)$ for two Banach spaces $E$ and $F$ the bounded finite rank operators from $E$ to $F$; as usual, we write ${\cal F}(E)$ as shorthand for ${\cal F}(E,E)$.
\begin{proposition} \label{facam}
Let $E$ and $F$ be Banach spaces such that $E^\ast$ and $F^\ast$ have the bounded approximation property. Suppose further that the following factorization property holds: 
\begin{quote}
there is $C \geq 0$ such that, for each $T \in {\cal F}(F)$, there are $S \in {\cal F}(F,E)$ and $R \in {\cal F}(E,F)$ with $\| R \| \| S \| \leq C \| T \|$ and $RS = T$.
\end{quote}
Then, for any index set $\mathbb I$, the following are equivalent:
\begin{items} 
\item $\ell^\infty(\mathbb{I},{\cal K}(E))$ is amenable;
\item $\ell^\infty(\mathbb{I},{\cal K}(E \oplus F))$ is amenable.
\end{items}
\end{proposition}
\begin{proof}
Let $\mathbb{I}$ be an index set, and let $\A := \ell^\infty(\mathbb{I},{\cal K}(E \oplus F))$. We wish to apply Theorem \ref{bgn}.
\par 
As both $E^\ast$ and $F^\ast$ have the bounded approximation property, $(E \oplus F)^\ast \cong E^\ast \oplus F^\ast$ has it, too. Consequently, ${\cal K}(E \oplus F)$ has a bounded approximate identity, and so does $\A$. 
\par 
For $i \in \mathbb{I}$, let $P_{1,i} \!: E \oplus F \to E$ and $P_{2,i} \!: E \oplus F \to F$ be the canonical projection; for $j=1,2$, set $P_j := ( P_{j,i})_{i \in \mathbb I} \in \ell^\infty(\mathbb{I},{\cal B}(E \oplus F)) \subset {\cal M}(\ell^\infty(\mathbb{I},{\cal K}(E \oplus F)))$. It follows that $P_1 \A P_1 \cong \ell^\infty(\mathbb{I},{\cal K}(E))$ and $P_2 \A P_2 \cong \ell^\infty(\mathbb{I},{\cal K}(F))$.
\par 
The restriction of $\Delta_\A$ to $P_2 \A P_1\Tensor P_1 \A P_2$ induces a quotient norm, say $| \cdot |$, on its range, which dominates the given norm $\| \cdot \|$. By our factorization hypothesis, the range of $\Delta_\A( P_2 \A P_1 \Tensor P_1 \A P_2)$ contains $\ell^\infty({\cal F}(F))$, and we have $| \cdot | \leq C \| \cdot \|$ on $\ell^\infty({\cal F}(F))$. As $F^\ast$ has the approximation property, so does $F$, and, in particular, ${\cal F}(F)$ is dense in ${\cal K}(F)$, as is $\ell^\infty({\cal F}(F))$ in $\ell^\infty({\cal K}(F))$. Every element of $\ell^\infty({\cal K}(F))$ is thus a limit---with respect to $\| \cdot \|$---of a sequence in $\ell^\infty({\cal F}(F))$. This sequence is a Cauchy sequence with respect to $\| \cdot \|$ and thus with respect to $| \cdot |$; consequently, it converges---with respect to $| \cdot |$---to an element in $\Delta_\A( P_2 \A P_1 \Tensor P_1 \A P_2)$. Since $\| \cdot \| \leq | \cdot |$, this limit with respect to $| \cdot |$ is the same limit as with respect to $\| \cdot \|$. So, $\Delta_\A$ maps $P_2 \A P_1 \Tensor P_1 \A P_2$ onto $P_2 \A P_2$.
\par
The hypotheses of Theorem \ref{bgn} are thus all satisfied, and the claim follows.
\end{proof}
\par 
We shall now look at Banach spaces for which the factorization hypothesis of Proposition \ref{facam} is satisfied.
\par 
Let $p \in [1,\infty]$ and $\lambda \geq 1$. A Banach space $E$ is called an \emph{${\cal L}^p_\lambda$-space} if, for any finite-dimensional subspace $X$ of $E$, there are $n \in \posints$ and an $n$-dimensional subspace $Y$ of $E$ containing $X$ such that $d(Y,\ell^p_n) \leq \lambda$, where $d$ is the Banach--Mazur distance (\cite[Definition 3.1]{LP}). We call $E$ simply an \emph{${\cal L}^p$-space} if it is an ${\cal L}^p_\lambda$-space for some $\lambda \geq 1$. All $L^p$-spaces, i.e., spaces of $p$-integrable functions on some measure space, are ${\cal L}^p$-spaces. The ${\cal L}^p$-spaces were introduced in \cite{LP} and studied further in \cite{LR}. We list some of their properties:
\begin{itemize}
\item Every ${\cal L}^p$-space is isomorphic to a subspace of an $L^p$-space (\cite[Theorem I(i)]{LR}). In particular, for $p \in (1,\infty)$, each ${\cal L}^p$-space is reflexive.
\item If $E$ is an ${\cal L}^p$-space, then $E^\ast$ is an ${\cal L}^{p'}$-space, where $p' \in [1,\infty]$ is conjugate to $p$, i.e., satisfies $\frac{1}{p} + \frac{1}{p'} = 1$ (\cite[Theorem III(a)]{LR}).
\item If $E$ is an ${\cal L}^p$-space, then there is a constant $\rho \geq 1$ such that, for each finite-dimensional subspace $X$ of $E$, there are $n \in \posints$, an $n$-dimensional subspace $Y$ of $E$ containing $X$ with $d(Y,\ell^p_n)$, and a projection $P$ onto $Y$ with $\| P \| \leq \rho$ (\cite[Theorem III(c)]{LR}). In particular, $E$ has the bounded approximation property.
\end{itemize}
\par 
In \cite{GJW}, it is mentioned without proof before \cite[Theorem 6.4]{GJW} that, given any two infinite-dimensional ${\cal L}^p$-spaces $E$ and $F$, every operator in ${\cal F}(F)$ factors through $E$ with both factors being compact. Since, for our purpose, we need control over the norms of those factors, we give a refinement of this observation with a detailed proof:
\begin{lemma} \label{Lpfac}
Let $p \in [1,\infty]$, and let $E$ and $F$ be ${\cal L}^p$-spaces with $\dim E = \infty$. Then there is $C \geq 0$ such that, for each $T \in {\cal F}(F)$, there are $S \in {\cal F}(F,E)$ and $R \in {\cal F}(E,F)$ with $\| R \| \| S \| \leq C \| T \|$ and $RS = T$.
\end{lemma}
\begin{proof}
Let $T \in {\cal F}(F)$, and set $X := TF$. Let $\lambda \geq 1$ be such that $E$ is a ${\cal L}^p_\lambda$-space. Then there are $n \in \posints$, a finite-dimensional subspace $Y$ of $F$ containing $X$, and a bijective linear map $\tau \!: Y \to \ell^p_n$ such that $\| \tau \| \| \tau^{-1} \| \leq \lambda$.
\par 
Let $\rho \geq 1$ be the constant for $E$ whose existence is guaranteed by \cite[Theorem III(c)]{LR}. Let $Z_0$ be an $n$-dimensional subspace of $E$. (Here, we require that $\dim E = \infty$.) Then there are $m \in \posints$, an $m$-dimensional subspace $Z$ of $E$ containing $Z_0$, a bijective map $\sigma \!: Z \to \ell^p_m$ with $\| \sigma \| \| \sigma^{-1} \| \leq \rho$, and a projection $P$ onto $Z$ with $\| P \| \leq \rho$; note that necessarily $m \geq n$. With $\iota \!: \ell^p_n \to \ell^p_m$ and $\pi \!: \ell^p_m \to \ell^p_n$ being the canonical embedding and projection, we see that $T=  \tau^{-1} \pi \sigma \sigma^{-1} \iota \tau T$. Set $S := \sigma^{-1} \iota \tau T$ and $R := \tau^{-1} \pi \sigma P$. Then $S=RT$ holds, and we have
\[
  \| R \| \| S \| \leq \| \tau^{-1} \| \| \pi \| \| \sigma \| \| P \| \| \sigma^{-1} \| \| \iota \| \| \tau \| \| T \|
  = \| \tau \| \| \tau^{-1} \| \| \sigma \| \| \sigma^{-1} \| \| P \| \| T \| \leq \lambda \rho^2 \| T \|.
\]
Hence, $C := \lambda \rho^2$ has the desired property.
\end{proof}
\begin{remark}
In \cite{DF}, A.\ Defant and K.\ Floret introduced the class of ${\cal L}^p_g$-spaces, which contains all ${\cal L}^p$-spaces, but is somewhat better behaved. For $p=1,\infty$, the ${\cal L}^p_g$-spaces are just the ${\cal L}^p$-spaces whereas, for $p \in (1,\infty)$, a space is an ${\cal L}^p_g$-space if and only if it is an ${\cal L}^p$-space \emph{or isomorphic to a Hilbert space} (\cite[23.3]{DF}). Therefore, Lemma \ref{Lpfac} does not hold true for ${\cal L}^p_g$-spaces if $p \neq 1,\infty$.
\end{remark}
\par 
Together, Proposition \ref{facam} and Lemma \ref{Lpfac} yield the following dichotomy theorem:
\begin{theorem} \label{thm3}
Let $p \in [1,\infty]$, and let $\mathbb I$ be an index set. Then one of the following assertions is true:
\begin{items}
\item $\ell^\infty(\mathbb{I},{\cal K}(E))$ is amenable for every infinite-dimensional ${\cal L}^p$-space $E$;
\item $\ell^\infty(\mathbb{I},{\cal K}(E))$ is not amenable for any infinite-dimensional ${\cal L}^p$-space $E$.
\end{items}
\end{theorem}
\begin{proof}
Suppose that (ii) is false, i.e, there is an infinite-dimensional ${\cal L}^p$-space $E$ such that $\ell^\infty(\mathbb{I},{\cal K}(E))$ is amenable. Let $F$ be any infinite-dimensional ${\cal L}^p$-space. Then Lemma \ref{Lpfac} and Proposition \ref{facam} imply that $\ell^\infty(\mathbb{I},{\cal K}(E \oplus F))$ is amenable. Interchanging the r\^oles of $E$ and $F$ and invoking Lemma \ref{Lpfac} and Proposition \ref{facam} again, then yields the amenability of $\ell^\infty(\mathbb{I},{\cal K}(F))$. This proves that (i) is true.
\end{proof}
\par
Combining Theorems \ref{thm2} and \ref{thm3}, we obtain:
\begin{corollary} \label{cor3}
Let $p \in (1,\infty]$, let $E = \ell^p$ if $p \in (1,\infty)$ and $E = c_0$ if $p = \infty$, and suppose that $\ell^\infty({\cal K}(E))$ is amenable. Then $\ell^\infty(\mathbb{I},{\cal K}(F))$ is amenable for every index set $\mathbb I$ and every infinite-dimensional ${\cal L}^p$-space $F$. In particular, ${\cal K}(F)$ is ultra-amenable for every infinite-dimensional ${\cal L}^p$-space $F$.
\end{corollary}
\begin{remark}
Let $p \in (1,\infty) \setminus \{ 2 \}$. Then $\ell^p(\ell^2)$ and $\ell^p \oplus \ell^2$ are not $L^p$-spaces, but still ${\cal L}^p$-spaces (\cite[Example 8.2]{LP}). Hence, if $\ell^\infty({\cal K}(\ell^p))$ is amenable, then so are $\ell^\infty({\cal K}(\ell^p(\ell^2)))$ and $\ell^\infty({\cal K}(\ell^p \oplus \ell^2))$. This is remarkable because ${\cal B}(\ell^p \oplus \ell^2)$ is known to be not amenable (\cite[Question 4]{Gro}).
\end{remark}
\section{A non-amenability criterion for $\ell^\infty({\cal K}(E \oplus F))$}
As we just observed, the amenability of $\ell^\infty({\cal K}(\ell^p))$ implies the amenability of $\ell^\infty({\cal K}(\ell^p \oplus \ell^2))$. In this section, we shall thus explore the amenability of $\ell^\infty({\cal K}(E \oplus F))$ for two Banach spaces $E$ and $F$.
\par 
Recall that an \emph{operator ideal} $\cal A$ is a rule that assigns to each pair $(E,F)$ of a Banach spaces a subspace ${\cal A}(E,F)$ of ${\cal B}(E,F)$ containing ${\cal F}(E,F)$ such that $RTS \in {\cal A}(X,Y)$ for any Banach spaces $X$ and $Y$, $T \in {\cal A}(E,F)$, $S \in {\cal B}(X,E)$, and $R \in {\cal B}(F,Y))$; if $E = F$, we convene again to simply write ${\cal A}(E)$. The seminal reference on operator ideals is \cite{Pie}. More recent treatments can be found in \cite{DF}, \cite{DJT}, or \cite{TJ}. We call $[{\cal A},\alpha]$---following \cite{DJT} in our notation---a \emph{Banach operator ideal} if, for each pair $(E,F)$ of Banach spaces, there is a norm $\alpha$ on ${\cal A}(E,F)$ turning it into a Banach space such that
\[
  \alpha(y \odot \phi) = \| y \| \| \phi \| \qquad (y \in F, \, \phi \in E^\ast),
\]
where $y \odot \phi \in {\cal F}(E,F)$ is the rank one operator corresponding to the elementary tensor $y \tensor \phi$, and
\[
  \alpha(RTS) \leq \| R \| \alpha(T) \| S \| \qquad (T \in {\cal A}(E,F), \, S \in {\cal B}(X,E), \, R \in {\cal B}(F,Y))
\]
for any Banach spaces $X$ and $Y$. Given a Banach operator ideal $[{\cal A},\alpha]$, its \emph{maximal hull} $[{\cal A}^{\max},\alpha^{\max}]$ is defined as follows. For two Banach spaces $E$ and $F$, let $\mathfrak{F}(E)$ denote the finite-dimensional subspaces of $E$, and let $\mathfrak{F}_c(F)$ stand for the closed subspaces of $F$ with finite co-dimension; for $X \in \mathfrak{F}(E)$ and $Y \in \mathfrak{F}_c(F)$, let $\iota_X \!: X \to E$ and $\pi_Y \!: F \to F/Y$ be the inclusion and quotient map, respectively. We define
\[
  \alpha^{\max}(T) := \sup \{ \alpha(\pi_Y T \iota_X ): X \in \mathfrak{F}(E), \, Y \in \mathfrak{F}_c(F) \} \in [0,\infty]
  \qquad (T \in {\cal B}(E,F))
\] 
and
\[
  {\cal A}^{\max}(E,F) := \left\{ T \in {\cal B}(E,F): \alpha^{\max}(T) < \infty \right\}.
\]
It is routinely checked that $[{\cal A}^{\max},\alpha^{\max}]$ is again a Banach operator ideal, and we call $[{\cal A},\alpha]$ \emph{maximal} if $[{\cal A}^{\max},\alpha^{\max}] = [{\cal A},\alpha]$.
\par 
An immediate consequence of \cite[17.5, Representation Theorem for Maximal Operator Ideals]{DF} is: if $[{\cal A},\alpha]$ is a maximal Banach operator ideal, then ${\cal A}(E,F^\ast)$ is a dual Banach space for any two Banach spaces $E$ and $F$ such that weak$^\ast$ topology of ${\cal A}(E,F^\ast)$ coincides with the weak$^\ast$ topology of ${\cal B}(E,F^\ast) = (E \Tensor F)^\ast$ on norm bounded subsets.
\par 
In particular, we have (compare \cite[17.21, Proposition]{DF}):
\begin{lemma} \label{maxid}
Let $[{\cal A},\alpha]$ be a maximal Banach operator ideal, let $E$ and $F$ be Banach spaces, and let $( T_i )_{i \in \mathbb I}$ be a bounded net in ${\cal A}(E,F^\ast)$ that converges to $T \in {\cal B}(E,F^\ast)$ with respect to the weak$^\ast$ topology of ${\cal B}(E,F^\ast)$. Then $T$ lies in ${\cal A}(E,F^\ast)$.
\end{lemma}
\begin{proposition} \label{opidnonam}
Let $E$ and $F$ be Banach spaces, let $\mathbb{I}$ be an index set, and suppose that there are $( T_i )_{i \in \mathbb I} \in \ell^\infty(\mathbb{I},{\cal K}(F,E))$ and an ultrafilter $\cal U$ over $\mathbb I$ such that $\text{\emph{weak}$^\ast$-}\lim_{i \in \cal U} T_i \notin {\cal K}(F,E^{\ast\ast})$, where the limit is with respect to the weak$^\ast$ topology of ${\cal B}(E,F^{\ast\ast})$. Suppose further that there is a maximal operator ideal $[{\cal A},\alpha]$ with the following properties:
\begin{alphitems}
\item ${\cal K}(E,F) \subset {\cal A}(E,F)$ such that the inclusion is continuous;
\item ${\cal A}(F,E^{\ast\ast}) \subset {\cal K}(F,E^{\ast\ast})$.
\end{alphitems}
Then $\ell^\infty(\mathbb{I}, {\cal K}(E \oplus F))$ is not amenable.
\end{proposition}
\begin{proof}
Let $P_1$ and $P_2$ be the projections in $\ell^\infty(\mathbb{I}, {\cal B}(E \oplus F))$ induced by the canonical projections onto $E$ and $F$, respectively (compare the proof of Proposition \ref{facam}). We wish to apply Proposition \ref{noamprop}.
\par 
Letting
\[
  I := \varcl{\ell^\infty(\mathbb{I},{\cal A}(E \oplus F)) \cap \ell^\infty(\mathbb{I},{\cal K}(E \oplus F))},
\]
with the closure taken in the norm topology of $\ell^\infty(\mathbb{I},{\cal K}(E \oplus F))$, defines a closed ideal of $\A := \ell^\infty(\mathbb{I},{\cal K}(E \oplus F))$.
\par 
From (a), it is obvious that $P_2 I P_1 = P_2 \A P_1$. 
\par
To see that $P_1 I P_2 \subsetneq P_1 \A P_2$, let $( T_i )_{i \in \mathbb I} \in \ell^\infty(\mathbb{I},{\cal K}(F,E))$ and an ultrafilter $\cal U$ over $\mathbb I$ be as specified in the hypotheses. We claim that $(T_i )_{i \in \mathbb I} \in P_1 \A P_2 \setminus P_1 I P_2$. Define
\[
  Q_{\cal U} \!:  \ell^\infty(\mathbb{I},{\cal K}(F,E)) \to {\cal B}(F,E), \quad
  ( S_i )_{i \in \mathbb I} \mapsto \text{weak$^\ast$-}\lim_{i \in \cal U} S_i.
\]
so that $T := Q_{\cal U}\left( (T_i )_{i \in \mathbb I} \right) \notin {\cal K}(F,E^{\ast\ast})$. Assume that $(T_i )_{i \in \mathbb I} \in P_1 I P_2$, and let $\epsilon >0$ be arbitrary. By the definition of $I$, there is thus $( R_i )_{i \in \mathbb I} \in 
\ell^\infty(\mathbb{I},{\cal A}(F,E)) \cap \ell^\infty(\mathbb{I},{\cal K}(F,E))$ such that $\sup_{i \in \mathbb I} \| R_i - T_i \| < \epsilon$. Since $Q_{\cal U}$ is a contraction, this means that $\| R - T \| < \epsilon$, where $R := Q_{\cal U}\left( (R_i )_{i \in \mathbb I} \right)$. By Lemma \ref{maxid}, $R \in {\cal A}(F,E^{\ast\ast})$ holds, so that $R \in {\cal K}(F,E^{\ast\ast})$ by (b). Since $\epsilon > 0$ was arbitrary, this means that $T \in {\cal K}(F,E^{\ast\ast})$, which contradicts our hypotheses.
\end{proof}
\par 
The hypotheses of Proposition \ref{opidnonam} appear to be technical and somewhat contrived, but as our next theorem shows, they do, in fact, occur naturally in certain situations:
\begin{theorem} \label{thm4}
The Banach algebra $\ell^\infty({\cal K}(E \oplus \ell^2))$ is not amenable for any of the following spaces $E$: $c_0$, $\ell^\infty$, and $\ell^1$.
\end{theorem}
\begin{proof}
We first consider the case $E = c_0$.
\par 
For $n \in \posints$, let $\pi_n \!: c_0 \to c_0$ denote the projection onto the first $n$ coordinates, and let $\iota \!: \ell^2 \to c_0 \hookrightarrow \ell^\infty$ be the natural inclusion. Then $( \pi_n \iota )_{n=1}^\infty \in \ell^\infty({\cal K}(\ell^2,c_0))$, and $\text{weak$^\ast$-}\lim_{n \in \cal U} \pi_n \iota = \iota \notin {\cal K}(\ell^2,\ell^\infty)$ for any free ultrafilter $\cal U$ over $\posints$.
\par 
Let $[\Pi_2,\pi_2]$ be the ideal of $2$-summing operators (see \cite[p.\ 31]{DJT}, for instance); then $[\Pi_2,\pi_2]$ is maximal (by \cite[6.16 Theorem]{DJT}). By \cite[3.7 Theorem]{DJT}, $\Pi_2(c_0,\ell^2) = {\cal B}(c_0,\ell^2)$ holds, so that Proposition \ref{opidnonam}(a) is satisfied. On the other hand, every $2$-summing operator is completely continuous (\cite[2.17 Theorem]{DJT}). Since $\ell^2$ is reflexive, this means that $\Pi_2(\ell^2,\ell^\infty) \subset {\cal K}(\ell^2,\ell^\infty)$, so that Proposition \ref{opidnonam}(b) holds as well.
\par 
The $E = \ell^\infty$ case has an almost identical proof.
\par 
Suppose now that $E = \ell^1$. We shall apply Proposition \ref{opidnonam} to the Banach algebra $\ell^\infty({\cal K}(\ell^2 \oplus \ell^1 ))$, which is isomorphic to $\ell^\infty({\cal K}(\ell^1 \oplus \ell^2))$. With $\iota \!: \ell^1 \to \ell^2$ being the canonical inclusion and $\pi_n \!: \ell^2 \to \ell^2$ for $n \in \posints$ denoting the projection onto the first $n$ coordinates, we have $\text{weak$^\ast$-}\lim_{n \in \cal U} \pi_n \iota = \iota \notin {\cal K}(\ell^1,\ell^2)$ for any free ultrafilter $\cal U$ over $\posints$---just as in the case $E = c_0$. Let Let $[\Pi_2^d,\pi_2^d]$ be the dual ideal of $[\Pi_2,\pi_2]$ (see \cite[p.\ 186]{DJT}), so that $[\Pi_2^d,\pi_2^d]$ is maximal by \cite[9.4 Corollary]{DJT}. Using the definition of $[\Pi_2^d,\pi_2^d]$ and arguing as in the $E= c_0$ case, we see that Proposition \ref{opidnonam}(a) and (b) are satisfied.
\end{proof}
\begin{remark}
Even though $\ell^\infty({\cal K}(\ell^2))$ is known not to be amenable, there seems to be no way---by means of Theorem \ref{bgn}, for instance---to conclude directly from its non-amenability that the Banach algebras considered in Theorem \ref{thm4} are not amenable.
\end{remark}
\par 
Having established the non-amenability of $\ell^\infty({\cal K}(\ell^p \oplus \ell^2))$ for $p =1,\infty$ with the help of Proposition \ref{opidnonam}, one might be tempted to try to extend this result to general $p \in [1,\infty] \setminus \{ 2 \}$ through the choice of a suitable maximal Banach operator ideal $[ {\cal A}, \alpha ]$. Alas, as we shall see now, this attempt is futile:
\begin{proposition} \label{lp2}
Let $p \in (1,\infty) \setminus \{ 2 \}$, let $\A := \ell^\infty({\cal K}(\ell^p \oplus \ell^2))$, and let $P_1, P_2 \in \ell^\infty({\cal B}(\ell^p \oplus \ell^2))$ be the projections induced by the canonical projections onto $\ell^p$ and $\ell^2$, respectively. Then there is no closed ideal $I$ of $\A$ with $P_2 I P_1 = P_2 \A P_1$ and $P_1 I P_2 \subsetneq P_1 \A P_2$.
\end{proposition}
\begin{proof}
We use the fact (\cite[p.\ 73]{LT}) that we have an isomorphism
\begin{equation} \label{iso}
  \ell^p \cong \text{$\ell^p$-}\bigoplus_{n=1}^\infty \ell^2_n.
\end{equation}
For $n \in \posints$, let $J_n \!: \ell^2_n \to \ell^p$ and $Q_n \!: \ell^p \to \ell^2_n$ the embedding of and projection onto the $n$-th summand in (\ref{iso}), respectively; note that those maps $J_n$ and $P_n$ are uniformly bounded. Furthermore, let $\iota_n \!: \ell^2_n \to \ell^2$ and $\pi_n \!: \ell^2 \to \ell^2_n$ be the canonical embedding and projection, respectively, for $n \in \posints$. 
\par 
Assume that there is a closed ideal $I$ of $\A$ with $P_2 I P_1 = P_2 \A P_1$. Note that, since $\A$ has a bounded approximate identity, $I$ is also a closed ideal of $\ell^\infty({\cal B}(\ell^p \oplus \ell^2))$, so that, in particular, $P_j I P_k \subset I$ for $j,k=1,2$.
\par 
Let $( T_n )_{n=1}^\infty \in \ell^\infty({\cal K}(\ell^2,\ell^p))$, and let $\epsilon > 0$ be arbitrary. For each $n \in \posints$, we can find $S_n \in {\cal F}(\ell^2,\ell^p)$ with $\| T_n - S_n \| \leq \epsilon$ as well as $N_n \in \posints$ such that $S_n \iota_{N_n} \pi_{N_n} = S_n$. Define $( U_n )_{n=1}^\infty \in \ell^\infty({\cal K}(\ell^2,\ell^p))$ and $( V_n )_{n=1}^\infty \in \ell^\infty({\cal K}(\ell^p))$ by letting
\[
  U_n := \iota_{N_n} Q_{N_n}  \quad\text{and}\quad V_n := J_{N_n}\pi_{N_n} \qquad (n \in \posints).
\]
so that
\[
  \begin{bmatrix} 0 & ( S_n )_{n=1}^\infty \\ 0 & 0 \end{bmatrix}
  = \begin{bmatrix} 0 & ( S_n )_{n=1}^\infty \\ 0 & 0 \end{bmatrix}
  \begin{bmatrix} 0 & 0 \\ ( U_n)_{n=1}^\infty & 0 \end{bmatrix}
  \begin{bmatrix} 0 & ( V_n )_{n=1}^\infty \\ 0 & 0 \end{bmatrix}.
\]
As $\left[ \begin{smallmatrix} 0 & 0 \\ ( U_n)_{n=1}^\infty & 0 \end{smallmatrix} \right] \in P_2 \A P_1 = P_2 I P_1 \subset I$ and $I$ is an ideal, it follows that $\left[ \begin{smallmatrix} 0 & ( S_n )_{n=1}^\infty \\ 0 & 0 \end{smallmatrix} \right] \in I$ as well. Since $\epsilon > 0$ was arbitrary, this entails that $( T_n )_{n=1}^\infty \in P_1 I P_2$, and since $( T_n )_{n=1}^\infty \in \ell^\infty({\cal K}(\ell^2,\ell^p))$ was arbitrary, this means that $P_1 I P_2 = P_1 \A P_2$.
\end{proof}
\begin{remarks}
\item Even though Proposition \ref{lp2} shows that the (still hypothetical) non-amenability of $\ell^\infty({\cal K}(\ell^p \oplus \ell^2))$ cannot be established in the same way as for ${\cal B}(\ell^p \oplus \ell^2)$, it naturally leads to the question if, for sufficiently nice Banach spaces $E$, the amenability of $\ell^\infty({\cal K}(E))$ forces ${\cal B}(E)$ to be amenable. Since ${\cal K}(E)^{\ast\ast} = {\cal B}(E)$ via trace duality for any reflexive Banach space with the approximation property, this question can, for such spaces, be put into a more general framework: If $\ell^\infty(\A)$ is amenable for some Banach algebra $\A$, does this imply that $\A^{\ast\ast}$, equipped with one of the Arens products (see \cite{Dal}), is amenable? Partial answers, which do not apply to the case where $\A = {\cal K}(\ell^p \oplus \ell^2)$, are given in \cite{CSR}. 
\item In \cite[Theorem 1]{CSR}, the following is claimed to be a consequence of \cite{GI}: For a unital $\cstar$-algebra $\A$, there are an index set $\mathbb I$, which can be chosen as $\posints$ if $\A^\ast$ is separable, and a an algebra homomorphism from $\ell^\infty(\mathbb{I},\A)$ onto $\A^{\ast\ast}$. An inspection of the proof of \cite[Theorem 1]{CSR} shows that the alleged algebra homomorphism is
\begin{equation} \label{bogus}
  \ell^\infty(\mathbb{I},\A) \to \A^{\ast\ast}, \quad ( a_i )_{i \in \mathbb I} \mapsto \text{weak$^\ast$-}\lim_{i \in \cal U} a_i
\end{equation}
for a suitable ultrafilter $\cal U$ over $\mathbb I$. Let $\A$ be the unitization of ${\cal K}(\ell^2)$, and let $\cal U$ be a free ultrafilter over $\posints$. Then we have 
\[
  \text{weak$^\ast$-}\lim_{n \in \cal U} \delta_1 \odot \delta_n = \text{weak$^\ast$-}\lim_{n \in \cal U} \delta_n \odot \delta_1 = 0 
\]
whereas
\[ 
  \text{weak$^\ast$-}\lim_{n \in \cal U} (\delta_1 \odot \delta_n)(\delta_n \odot \delta_1) = \delta_1 \odot \delta_1 \neq 0,
\]
which means that (\ref{bogus}) cannot be multiplicative, contrary to what is claimed in \cite{CSR}.
\end{remarks}
\section{A left ideal in $({\cal K}(\ell^p))_{\cal U}$ without bounded right approximate identity}
If $\A$ is an amenable Banach algebra, then a closed ideal of $\A$ has a bounded approximate identity if and only if it is weakly complemented. By finding a closed, weakly complemented ideal of $\A$ that lacks a bounded approximate identity, one can thus show that $\A$ is not amenable, as was done in \cite{DGH} in the case of the measure algebra of a non-discrete, locally compact group. More generally, a closed left ideal of $\A$ has a bounded right approximate identity if and only if it is weakly complemented (see \cite[Lemma 2.3.6]{LoA}, for instance).
\par
In this section, for $p \in (1,2)$ and a free ultrafilter $\cal U$ over $\posints$, we shall exhibit a closed left ideal of $({\cal K}(\ell^p))_{\cal U}$ that lacks a right approximate identity (bounded or not) and present some, albeit circumstantial, evidence for it being weakly complemented.
\par 
For $p \in (1,2)$, let $\iota \!: \ell^p \to \ell^2$ be the natural inclusion map, and note that the adjoint $\iota^\ast \!: \ell^2 \to \ell^{p'}$ is the canonical inclusion of $\ell^2$ in $\ell^{p'}$. Let $\cal U$ be a free ultrafilter over $\posints$, and define
\begin{equation} \label{leftid}
  L_2 := \varcl{\{ ( T_n \iota )_{\cal U} : ( T_n )_{\cal U} \in ({\cal B}(\ell^2,\ell^p))_{\cal U} \}}^{({\cal B}(\ell^p))_{\cal U}}.
\end{equation}
Obviously, $L_2$ is a closed left ideal of $({\cal B}(\ell^p))_{\cal U}$. By Pitt's theorem (\cite[Proposition 2.c.3]{LT}), ${\cal B}(\ell^2,\ell^p) = {\cal K}(\ell^2,\ell^p)$ holds, so that $L_2$ is, in fact, a closed ideal of $({\cal K}(\ell^p))_{\cal U}$.
\par 
Recall that, for $\epsilon > 0$, a $(1+\epsilon)$-isometry from a Banach space $E$ into a Banach space $F$, is a linear map $T \!: E \to F$ satisfying
\[
  (1 -\epsilon)\| x \| \leq \| Tx \| \leq (1+\epsilon) \| x \| \qquad (x \in E).
\]
By \cite[19.1 Dvoretzky's Theorem]{DJT}, there is, for each infinite-dimensional Banach space $E$, for each $n \in \posints$, and for each $\epsilon > 0$, a $(1+\epsilon)$-isometry from $\ell^2_n$ into $E$. We shall use this theorem to obtain particular elements of $L_2$. For each $n \in \posints$, let $\pi_n \!: \ell^2 \to \ell^2_n$ denote the canonical projection onto the first $n$ coordinates, and let, for each $n \in \posints$, $\tau_n \!: \ell^2_n \to \ell^p$ be a $\left(1 + \frac{1}{n} \right)$-isometry, which exists by Dvoretzky's theorem. Then $(\tau_n \pi_n \iota )_{\cal U}$ lies in $L_2$.
\par 
The following is our technical main result in this section:
\begin{lemma} \label{dvolem}
Let $p \in (1,2)$, let $\cal U$ be a free ultrafilter over $\posints$, and let $L_2$ and $(\tau_n \pi_n \iota )_{\cal U}$ be defined as above. Then we have
\[
  \| (\tau_n \pi_n \iota )_{\cal U} (T_n)_{\cal U} - (\tau_n \pi_n \iota )_{\cal U} \| \geq 1
  \qquad ((T_n)_{\cal U} \in L_2).
\]
\end{lemma}
\begin{proof}
Assume towards a contradiction that there are $\theta \in [0,1)$ and $( T_n )_{\cal U} \in ({\cal B}(\ell^2,\ell^p))_{\cal U}$ such that
\begin{equation} \label{dvoineq}
  \theta > \| (\tau_n \pi_n \iota )_{\cal U} (T_n \iota)_{\cal U} - (\tau_n \pi_n \iota )_{\cal U} \|
  = \lim_{n \in \cal U} \| \tau_n \pi_n \iota T_n \iota - \tau_n \pi_n \iota \| 
  = \lim_{n \in \cal U} \| \pi_n \iota T_n \iota - \pi_n \iota \|,
\end{equation}
where the last equality is due to the fact that $\tau_n$ is a $\left(1 + \frac{1}{n} \right)$-isometry for each $n \in \posints$.
\par 
Let $T := \text{weak-}\lim_{n \in \cal U} T_n \in {\cal K}(\ell^2,\ell^p)$. (The limit exists by \cite[Propositon 1.45]{Hei} because ${\cal K}(\ell^2,\ell^p)$ is reflexive, so that its closed unit ball is weakly compact.) Let $x \in \ell^p$ and $\xi \in \ell^2$. Then we have:
\[
  \begin{split}
  | \langle \xi, (\iota T \iota - \iota)(x) \rangle | 
  & = | \langle \iota^\ast(\xi), (T\iota)(x) \rangle - \langle \xi, \iota(x) \rangle | \\
  & = \lim_{n \in \cal U} | \langle \iota^\ast(\xi), (T_n \iota)(x) \rangle - \langle \xi, \iota(x) \rangle | \\
  & = \lim_{n \in \cal U} | \langle \xi, (\iota T_n \iota - \iota)(x) \rangle | \\
  & = \lim_{n \in \cal U} | \langle \pi_n^\ast(\xi), (\iota T_n \iota - \iota)(x) \rangle |, 
    \qquad\text{because $\lim_{n \to \infty} \| \pi_n^\ast(\xi) - \xi \|_2 = 0$}, \\
  & = \lim_{n \in \cal U} | \langle \xi, (\pi_n\iota T_n \iota - \pi_n\iota)(x) \rangle | \\
  & \leq \lim_{n \in \cal U} \| \pi_n\iota T_n \iota - \pi_n\iota \| \| x \| \| \xi \| \\
  & < \theta \| x \| \| \xi \|, \qquad\text{by (\ref{dvoineq})}.
  \end{split}
\]
It follows that
\begin{equation} \label{dvoineq2}
  \| \iota T \iota - \iota \| \leq \theta.
\end{equation}
\par 
Since $T$ is compact, so is $\iota T \iota$. Hence, there is a strictly increasing sequence $( n_k)_{k=1}^\infty$ in $\posints$ such that $((\iota T \iota)(\delta_{n_k}))_{k=1}^\infty$ is norm convergent in $\ell^2$ with limit $\eta$, say. It follows that
\begin{equation} \label{dvoeq}
  \lim_{k \to \infty} \langle \delta_{n_k}, (\iota T \iota)(\delta_{n_k}) \rangle = 
  \lim_{k \to \infty} \langle \delta_{n_k},\eta \rangle = 0.
\end{equation}
Together, (\ref{dvoineq2}) and (\ref{dvoeq}) yield
\begin{multline*}
  1 = \lim_{k \to \infty} \langle \delta_{n_k}, \iota(\delta_{n_k}) \rangle \\
    = \lim_{k \to \infty} | \langle \delta_{n_k}, (\iota T \iota)(\delta_{n_k}) \rangle -
                            \langle \delta_{n_k}, \iota(\delta_{n_k}) \rangle |
    = \lim_{k \to \infty} | \langle \delta_{n_k}, (\iota T \iota - \iota)(\delta_{n_k}) \rangle | \leq \theta,
\end{multline*}
which is impossible because $\theta \in [0,1)$.
\end{proof}
\par 
The following is now immediate:
\begin{proposition} \label{norai}
Let $p \in (1,2)$, let $\cal U$ be a free ultrafilter over $\posints$, and let $L_2$ be the closed left ideal of $({\cal K}(\ell^p))_{\cal U}$ defined in \emph{(\ref{leftid})}. Then $L_2$ does not have a right approximate identity.
\end{proposition}
\begin{remark}
Both Lemma \ref{dvolem} and Proposition \ref{norai} remain true in the slightly more general situation where $\cal U$ is a countably incomplete ultrafilter over an arbitrary index set.
\end{remark}
\par
\emph{If} we could establish that $L_2$ is weakly complemented, i.e., has a complemented annihilator in $({\cal K}(\ell^p))_{\cal U}^\ast$, then we know that $({\cal K}(\ell^p))_{\cal U}$---and thus, by Theorem \ref{thm1}, ${\cal B}(\ell^p)$---cannot be amenable. Unfortunately, such a proof eludes us, mostly due to the lack of a suitable description of $({\cal K}(\ell^p))_{\cal U}^\ast$. Nevertheless, we are able to show that the annihilator of $L_2$ in a certain closed subspace of $({\cal K}(\ell^p))_{\cal U}^\ast$ is indeed complemented.
\par 
We achieve this as a by-product of a general complementation result for ultrapowers of vector valued $\ell^p$-spaces.
\par 
Given a set $S$ and an ultrafilter $\cal U$ over some index set $\mathbb{I}$, we use $\langle S \rangle_{\cal U}$ for the corresponding \emph{set theoretic ultrapower} (see \cite{Hei} for the definition). For $( s_i )_{i \in \mathbb I} \in S^\mathbb{I}$, we denote its image in $\langle S \rangle_{\cal U}$ by $\langle s_i \rangle_{\cal U}$.
\par 
The following lemma relates, for $p \in [1,\infty)$, a Banach space $E$, and an ultrafilter $\cal U$ (over an arbitrary index set), the spaces $\ell^p(\langle \posints \rangle_{\cal U}, (E)_{\cal U})$ and $(\ell^p(E))_{\cal U}$. We identify the finitely supported functions in $\ell^p(\langle \posints \rangle_{\cal U}, (E)_{\cal U})$ with the algebraic tensor product $c_{00}(\langle \posints \rangle_{\cal U}) \tensor (E)_{\cal U}$, where $c_{00}(\langle \posints \rangle_{\cal U})$ are the finitely supported functions from $\langle \posints \rangle_{\cal U}$ into $\comps$.
\begin{lemma} \label{subspace}
Let $p \in [1,\infty)$, let $E$ be a Banach space, and let $\cal U$ be an ultrafilter. Then there is a unique isometry $J_p \!: \ell^p(\langle \posints \rangle_{\cal U}, (E)_{\cal U}) \to (\ell^p(E))_{\cal U}$ given by
\begin{equation} \label{embed}
  J_p \left( \delta_{\langle n_i \rangle_{\cal U}} \tensor ( x_i)_{\cal U} \right) = ( \delta_{n_i} \tensor x_i )_{\cal U}
  \qquad \left(\delta_{\langle n_i \rangle_{\cal U}} \in \langle \posints \rangle_{\cal U}, \, (x_i)_{\cal U} \in (E)_{\cal U} \right)
\end{equation}
\end{lemma}
\begin{proof}
It is routinely checked that (\ref{embed}) defines an isometry from $c_{00}(\langle \posints \rangle_{\cal U}) \tensor (E)_{\cal U}$ into $(\ell^p(E))_{\cal U}$, which then extends to all of $\ell^p(\langle \posints \rangle_{\cal U}, (E)_{\cal U})$ by continuity.
\end{proof}
\par 
Lemma \ref{subspace} enables us to canonically identify $\ell^p( \langle \posints \rangle_{\cal U},(E)_{\cal U})$ with a closed subspace of $(\ell^p(E))_{\cal U}$. 
\par 
Given a Banach space $E$ and an ultrafilter $\cal U$, there is a canonical duality between $(E)_{\cal U}$ and $(E^\ast)_{\cal U}$, which induces an isometric embedding of $(E^\ast)_{\cal U}$ into $(E)_{\cal U}^\ast$; for countably incomplete $\cal U$, this embedding is an isomorphism if and only if $(E)_{\cal U}$ is reflexive (\cite[Proposition 7.1]{Hei}). Recall that $E$ is called \emph{superreflexive} if every Banach space that can be finite represented in $E$ is reflexive; equivalently, $E$ is superreflexive if and only if $(E)_{\cal U}$ is reflexive for each ultrafilter $\cal U$ (\cite[Proposition 6.4]{Hei}). Also, if $E$ is superreflexive and $p \in (1,\infty)$, then $\ell^p(E)$ is also superreflexive (\cite[Proposition 4]{Daw0}). All this guarantees that the map $\Pi_p$ in the following proposition is well defined.
\begin{proposition} \label{kernel}
Let $p \in (1,\infty)$, let $E$ be a superreflexive Banach space, let $\cal U$ be an ultrafilter, and let $J_p \!: \!: \ell^p(\langle \posints \rangle_{\cal U}, (E)_{\cal U}) \to (\ell^p(E))_{\cal U}$ and $J_{p'} \!: \!: \ell^{p'}(\langle \posints \rangle_{\cal U}, (E^\ast)_{\cal U}) \to (\ell^{p'}(E^\ast))_{\cal U}$ be as in Lemma \emph{\ref{subspace}}. Then $\Pi_p := J_p J^\ast_{p'}$ is a norm one projection onto $\ell^p( \langle \posints \rangle_{\cal U},(E)_{\cal U})$. Moreover, we have for any $q \in (p,\infty]$ that
\begin{equation} \label{kernelq}
  \ker \Pi_p = \left\{ (x_i)_{\cal U} \in (\ell^p(E))_{\cal U} : \lim_{i \in \cal U} \| x_i \|_{\ell^q(E)} = 0 \right\}.
\end{equation}
\end{proposition}
\begin{proof}
It is easy to see that $\Pi_p$ is indeed a norm one projection onto $\ell^p( \langle \posints \rangle_{\cal U},(E)_{\cal U})$. 
\par 
Let $\mathbb I$ be the index set over which $\cal U$ is defined. For $( n_i )_{i \in \mathbb I} \in \posints^\mathbb{I}$, let $P_{n_i} \!: \ell^p(E) \to E$ denote the projection onto the $n_i$-th coordinate. From the definition of $\Pi_p$, it is clear that $( x_i )_{\cal U} \in (\ell^p(E))_{\cal U}$ belongs to $\ker \Pi_p$ if and only if $\lim_{i \in \cal U} \| P_{n_i} x_i \|_E = 0$ for any
$( n_i )_{i \in \mathbb I} \in \posints^\mathbb{I}$. It follows that
\[
  \ker \Pi_p \supset \left\{ (x_i)_{\cal U} \in (\ell^p(E))_{\cal U} : \lim_{i \in \cal U} \| x_i \|_{\ell^\infty(E)} = 0 \right\}
\]
For the converse inclusion, let $(x_i)_{\cal U} \in (\ell^p(E))_{\cal U}$ be such that $\lim_{i \in \cal U} \| x_i \|_{\ell^\infty(E)} =: \delta > 0$. Let $U \in \cal U$ be such that $\| x_i \|_{\ell^\infty(E)} = \sup_{n \in \posints} \| P_n x_i \| > \frac{\delta}{2}$ for each $i \in U$. For each $i \in U$, choose $n_i \in \posints$ such that $\| P_{n_i} x_i \|_{\ell^\infty(E)}  > \frac{\delta}{2}$. It follows that $\lim_{i \in \cal U} \| P_{n_i} x_i \|_{\ell^\infty(E)} \geq \frac{\delta}{2} > 0$, so that $(x_i)_{\cal U} \in (\ell^p(E))_{\cal U} \notin \ker \Pi_p$. All in all, we have
\begin{equation} \label{kernelinfty}
  \ker \Pi_p  = \left\{ (x_i)_{\cal U} \in (\ell^p(E))_{\cal U} : \lim_{i \in \cal U} \| x_i \|_{\ell^\infty(E)} = 0 \right\}
\end{equation}
\par 
Let $q \in (p,\infty)$. In view of (\ref{kernelinfty}), it is clear that 
\[
  \ker \Pi_p \supset \left\{ (x_i)_{\cal U} \in (\ell^p(E))_{\cal U} : \lim_{i \in \cal U} \| x_i \|_{\ell^q(E)} = 0 \right\}.
\]
For the converse inclusion, note that, for any $x = ( x_n )_{n=1}^\infty \in \ell^p(E)$, we have
\[
  \| x \|_{\ell^q(E)}^q = \sum_{n=1}^\infty \| x_n \|^q = 
  \sum_{n=1}^\infty \| x_n \|^{q-p} \| x_n \|^p 
  \leq \| x \|_{\ell^\infty(E)}^{q-p} \sum_{n=1}^\infty \| x_n \|^p
  \leq \| x \|_{\ell^\infty(E)}^{q-p} \| x \|^p_{\ell^p(E)}.
\]
Consequently, if $( x_i )_{\cal U} \in \ker \Pi_p$, i.e., $\lim_{i \in \cal U} \| x_i \|_{\ell^\infty(E)} = 0$, then
$\lim_{i \in \cal U} \| x_i \|_{\ell^q(E)} = 0$ holds as well. This proves (\ref{kernelq}).
\end{proof}
\par 
Let $p \in (1,2)$, and let $\cal U$ be a free ultrafilter over $\posints$. We can canonically represent $({\cal B}(\ell^p))_{\cal U}$ on $(\ell^p)_{\cal U}$ by letting
\[
  ( T_n )_{\cal U} ( x_n )_{\cal U} = (T_n x_n)_{\cal U} \qquad ( ( T_n )_{\cal U} \in ({\cal B}(\ell^p))_{\cal U}, \,
  ( x_n )_{\cal U} \in (\ell^p)_{\cal U}).
\]
Clearly, $( T_n )_{\cal U} ( x_n )_{\cal U} = 0$ holds for all $( T_n )_{\cal U} \in L_2$ if and only if $\lim_{n \in \cal U} \| x_n \|_2 = 0$, i.e., $( x_n )_{\cal U} \in \ker \Pi_p$ by Proposition \ref{kernel}. 
\par 
The Banach space ${\cal B}((\ell^p)_{\cal U})$ has the canonical predual $(\ell^p)_{\cal U} \Tensor (\ell^{p'})_{\cal U}$, which, by \cite[Proposition 4.7]{Daw}, embeds isometrically into $(\ell^p \Tensor \ell^{p'})_{\cal U} = ({\cal K}(\ell^p)^\ast)_{\cal U}$ and thus into $({\cal K}(\ell^p))_{\cal U}^\ast$ (see \cite[p.\ 87]{Hei}). It therefore makes sense to speak of the annihilator of $L_2$ in $(\ell^p)_{\cal U} \Tensor (\ell^{p'})_{\cal U}$. 
\par 
In view of the foregoing we have:
\begin{corollary} \label{complcor}
Let $p \in (1,2)$, and let $\cal U$ be a free ultrafilter. Then the annihilator of $L_2$ in $(\ell^p)_{\cal U} \Tensor (\ell^{p'})_{\cal U}$ is its complemented subspace $\ker \Pi_p  \Tensor (\ell^{p'})_{\cal U}$, where $\Pi_p$ is the canonical projection from $(\ell^p)_{\cal U}$ onto $\ell^p(\langle \posints \rangle_{\cal U})$.
\end{corollary}
\begin{remark}
It would be interesting to know whether the annihilator of $L_2$ in $({\cal K}(\ell^p)^\ast)_{\cal U}$ is complemented: as $({\cal K}(\ell^p))_{\cal U}^\ast$ can be finitely represented in $({\cal K}(\ell^p)^\ast)_{\cal U}$ (\cite[Theorem 7.3]{Hei}), this would further support our belief that $L_2$ is weakly complemented.
\end{remark}
\renewcommand{\baselinestretch}{1.0}
\renewcommand{\baselinestretch}{1.2}
\vfill
\begin{tabbing} 
\textit{Second author's address}: \= Department of Mathematical and Statistical Sciences \kill
\textit{First author's address}:  \> Department of Pure Mathematics \\
                                  \> University of Leeds \\
                                  \> Leeds, LS2 9JT \\
                                  \> United Kingdom \\[\medskipamount]
\textit{E-mail}:                  \> \texttt{matt.daws@cantab.net} \\[\bigskipamount]  
\textit{Second author's address}: \> Department of Mathematical and Statistical Sciences \\
                                  \> University of Alberta \\
                                  \> Edmonton, Alberta \\
                                  \> Canada T6G 2G1 \\[\medskipamount]
\textit{E-mail}:                  \> \texttt{vrunde@ualberta.ca}
\end{tabbing}
\end{document}